# Choosing the optimal multi-point iterative method for the Colebrook flow friction equation – Numerical validation

**Pavel Praks [1,*] and Dejan Brkić [1,*]**

[1] European Commission, DG Joint Research Centre (JRC), Directorate C: Energy, Transport and Climate, Unit C3: Energy Security, Distribution and Markets, Via Enrico Fermi 2749, 21027 Ispra (VA), Italy and IT4Innovations National Supercomputing Center, VŠB - Technical University Ostrava, 17. listopadu 2172/15, 708 00 Ostrava, Czech Republic; ORCID id: https://orcid.org/0000-0002-3913-7800, Pavel.Praks@ec.europa.eu, Pavel.Praks@vsb.cz (P.P.)

[2] European Commission, DG Joint Research Centre (JRC), Directorate C: Energy, Transport and Climate, Unit C3: Energy Security, Distribution and Markets, Via Enrico Fermi 2749, 21027 Ispra (VA), Italy; ORCID id: https://orcid.org/0000-0002-2502-0601, dejanbrkic0611@gmail.com, dejanrgf@tesla.rcub.bg.ac.rs (D.B.)

* Both authors contributed equally to this study

**Abstract:** The Colebrook equation $\zeta$ is implicitly given in respect to the unknown flow friction factor $\lambda$; $\lambda = \zeta(Re, \varepsilon^*, \lambda)$ which cannot be expressed explicitly in exact way without simplifications and use of approximate calculus. Common approach to solve it is through the Newton-Raphson iterative procedure or through the fixed-point iterative procedure. Both requires in some case even eight iterations. On the other hand numerous more powerful iterative methods such as three-or two-point methods, etc. are available. The purpose is to choose optimal iterative method in order to solve the implicit Colebrook equation for flow friction accurately using the least possible number of iterations. The methods are thoroughly tested and those which require the least possible number of iterations to reach the accurate solution are identified. The most powerful three-point methods require in worst case only two iterations to reach final solution. The recommended representatives are Sharma-Guha-Gupta, Sharma-Sharma, Sharma-Arora, Džunić-Petković-Petković; Bi-Ren-Wu, Chun-Neta based on Kung-Traub, Neta, and Jain method based on Steffensen scheme. The recommended iterative methods can reach the final accurate solution with the least possible number of iterations. The approach is hybrid between iterative procedure and one-step explicit approximations and can be used in engineering design for initial rough, but also for final fine calculations.

**Keywords:** Colebrook equation; Colebrook-White; iterative methods; three-point methods; turbulent flow; hydraulic resistances; pipes; explicit approximations;

---

## 1. Introduction

The Colebrook function $\zeta$ is to date the most used relation in engineering practice for evaluation of flow friction $\lambda$ in pipes [1]. It is given in implicit form $\lambda = \zeta(Re, \varepsilon^*, \lambda)$, Equation (1):

$$x = -2 \cdot log_{10}\left(\frac{2.51 \cdot x}{Re} + \frac{\varepsilon^*}{3.71}\right) \quad (1)$$

In Equation (1), $x = \frac{1}{\sqrt{\lambda}}$ is introduced because of linearization of the unknown flow friction factor $\lambda$, $Re$ is the Reynolds number and $\varepsilon^*$is the relative roughness of inner pipe surface (all quantities are dimensionless). Because the Colebrook function $\zeta$ in the examined iterative methods sometimes needs to be evaluated in two or three points, $y$ and $z$ is used in the same meaning as $x$, where they are dimensionless parameters that depend on friction factor. Practical domain of applicability is for the Reynolds number, $4000 < Re < 10^8$ and for the relative roughness of inner pipe surface, $0 < \varepsilon^* < 0.05$.

The Colebrook equation is empirical [2] and hence its accuracy can be disputed (experiment of Nikuradse or recent experiments from Oregon or Princeton research groups [3]), but anyway it is widely accepted in engineering practice. It is based on experiment conducted by Colebrook and White with flow of air through set of pipes with different inner roughness [2]. Turbulent part of the Moody diagram [5,6] is based on the Colebrook equation. The exact solution to this equation does not exist with the exception of those through the Lambert $W$-function [6-8], but anyway further the Lambert $W$-function can be evaluated only approximately [9-11]. Many different explicit approximations of the Colebrook equation exist in the form $\lambda = \varsigma(Re, \varepsilon^{*})$, and they are with different degree of accuracy and complexity [12-14]. Many of approximations are based on internal iterative cycles [15-17] and therefore it is better to use more accurate iterative procedures if they require only few iterative steps. Calculation of flow through complex networks of pipes such as for water or gas distribution requires multiple evaluations of flow friction factor [18-23]. In general the less number of iterations required, the solution is more efficient with decreased burden for computers [24-28] (recent approach based on Padé polynomials shows how the computational burden can be minimized with the same number of used iterations [29]).

The most used iterative methods for solving the Colebrook equation are the Newton-Raphson method or its simplified version, the fixed point methods [30,31]. On the other hand, various iterative procedures for solving a single non-linear equations are available and here are tested in total 23 different methods: i) Fixed-point [30,31]; ii) Newton-Raphson [32,33]; Hansen-Patrick group of method [34] such as: iii) Halley [35], iv) Euler-Chebyshev [36,37], and v) Basto-Semiao-Calheiros method [38]; vi) Super Halley [39]; vii) Ostrowski method of King family [36,37,40]; viii) Kung-Traub [41]; ix) Maheshwari [42]; x) Hermite interpolation [43] based on Jarratt method [44]; xi) Khattri and Babajee [45]; xii) Murakami [46]; xiii) Neta [47]; xiv) Chun-Neta [48] based on Kung-Traub method [41]; xv) Wang-Liu [49] based on Hermit interpolation [43]; xvi) Bi-Ren-Wu [50]; xvii) Jain [51] based on Steffensen method [52-55]; xviii) Sharma-Arora [56]; xix) Džunić-Petković-Petković [57-59]; xx) Neta-Johnson [60] based on Jarrat method [44]; xxi) Cordero et al. [54.55]; xxii) Sharma-Sharma [61]; and xxiii) Sharma-Guha-Gupta method [62]. All methods [63] are thoroughly tested within the domain of applicability of the Colebrook equation.

Details about preparation of data for analysis of the iterative methods are given in Section 2, while the methods are shown in Section 3 along with numerical examples. Discussion and analysis are given in Section 4 while concluding remarks and recommendations in Section 5.

## 2. Preparation of Data for Analysis

Calculation of complex networks of pipes for distribution of water or natural gas is computationally demanded task where flow friction needs to be evaluated even few million times which can cause significant burden for computers. Moreover, a fast but reliable approximation of pipeline hydraulic is needed also for probabilistic risk assessment of pipeline networks. A large number of network simulations of random component failures and their combinations must be automatically evaluated and statistically analyzed in this case [20,22]. These tasks can cause significant burden for computing which even powerful computer resources cannot easily deal with. In general logarithmic functions as well non-integer power terms are extremely demanded and therefore their use should be discouraged [24-27]. Both are used frequently in explicit approximations, but non-integer power terms very rarely in iterative procedures. So the main goal of this paper is to use optimal iterative solution for the Colebrook equation which leads to the accurate solution within the domain of applicability of the Colebrook equation after the least possible number of iterations. The Colebrook equation is based on logarithmic law and therefore it is hardly possible that the logarithmic function can be eliminated entirely from the iterative procedure [29].

To maintain compatibility in this analysis with the study of Brkić [13] that evaluates accuracy of explicit approximations of the Colebrook equation, all iterative methods are tested in 740 uniformly distributed points over the domain of applicability of the Colebrook equation. Examples from five points of 740 in total are shown (five triplets $\{Re, \varepsilon^{*}\} \rightarrow \{\lambda\}$), of which three are randomly selected and two are from the most problematic zones which needs slightly more number of iterations to reach the demanded level of accuracy: 1) $Re = 3.78 \cdot 10^6$, $\varepsilon^* = 0.00854$; 2) $Re = 6.23 \cdot 10^4$, $\varepsilon^* = 0.012$; 3) $Re = 1.18 \cdot 10^7$, $\varepsilon^* = 0.032$; 4) $Re = 5.74 \cdot 10^7$, $\varepsilon^* = 0.0008$; and 5) $Re = 8.31 \cdot 10^3$, $\varepsilon^* = 0.024$. The corresponding final solutions are: 1) $x = 5.274511499$; 2) $x = 4.928634498$; 3) $x = 4.128359435$; 4) $x = 7.331277467$; and 5) $x = 4.22204103$. Examples 2 and 5 are those from the more problematic zones of the domain which in general require more iterative steps to reach the desired accuracy.

Initial starting point for the all examined iterative procedures is set to; $x_0 = 7.273626085$. This value is chosen after numerous tests over the domain and it is from the problematic zone which needs in general more iterative cycles to reach the final solution. The fixed initial starting point for the testing is more suitable option because the iterative procedures in that way can be compared in a better way. Although the starting point $x_0$ for the iterative procedure can be chosen using different formulas [29], numerous tests showed that any fixed value within the domain of applicability of the Colebrook equation leads to the final accurate solution without significant variation of the required number of iterations. From the examined literature, cases when the iterative procedure diverges, oscillates or converges outside the domain of applicability of the Colebrook equation are not reported.

For the purpose of all examined iterative methods, the Colebrook equation should be in the appropriate form as in Equation (2); for point $x$; $x \to F(x) = x - \zeta(x)$, where $\zeta$ is the functional symbol for the Colebrook equation. Sometimes it need to be evaluated in additional points $y \to F(y) = y - \zeta(y)$ and $z \to F(z) = z - \zeta(z)$. The first $F'$ and the second derivatives $F''$ are also needed for some method but in most cases only in point $x$. Additional symbols specific for the certain method is explained in Section 3.

$$\left.\begin{aligned} x \to F(x) &= x + 2 \cdot log_{10}\left(\frac{2.51 \cdot x}{Re} + \frac{\varepsilon^*}{3.71}\right) \\ y \to F(y) &= y + 2 \cdot log_{10}\left(\frac{2.51 \cdot y}{Re} + \frac{\varepsilon^*}{3.71}\right) \\ z \to F(z) &= z + 2 \cdot log_{10}\left(\frac{2.51 \cdot z}{Re} + \frac{\varepsilon^*}{3.71}\right) \\ x \to F'(x) &= \frac{5.02}{2.3026 \cdot Re \cdot \left(\frac{10 \cdot \varepsilon^*}{37.1} + \frac{2.51 \cdot x}{Re}\right)} + 1 \\ x \to F''(x) &= -\frac{12.6}{2.3026 \cdot Re^2 \cdot \left(\frac{10 \cdot \varepsilon^*}{37.1} + \frac{2.51 \cdot x}{Re}\right)^2} \end{aligned}\right\} \tag{2}$$

In Equation (2), $\ln(10) \approx 2.3026$.

In Section 3, the presented iterative methods are listed in general from the simplest to the more complex [63]; 3.1 One log-call per iteration: 3.1.1) Fixed-point [30,31]; 3.1.2) Newton-Raphson [32,33]; 3.1.3) Hansen-Patrick [34]: 3.1.3a) Halley [35], 3.1.3b) Euler-Chebyshev [36,37], 3.1.3c) Basto-Semiao-Calheiros method [38]; 3.1.4) Super Halley [39]; 3.1.5) Murakami [46]; 3.2 Two log-calls per iteration_(two-point methods): 3.2.1) Ostrowski method of King family [36,37,40]; 3.2.2) Kung-Traub [41]; 3.2.3) Maheshwari [42]; 3.2.4) Khattri and Babajee [45]; 3.2.5) Hermite interpolation [43] based on Jarratt method [44]; 3.2.6) Wang-Liu method [49] based on Hermit interpolation [43]; and finally 3.3 Three-log-calls per iteration (three-point methods): 3.3.1) Neta [47]; 3.3.2) Chun-Neta [48] based on Kung-Traub [41]; 3.3.3) Džunić-Petković-Petković method [57-59]; 3.3.4) Neta-Johnson [60] based on Jarrat method [44]; 3.3.5) Jain [51] based on Steffensen method [52-55]; 3.3.6) Bi-Ren-Wu [50]; 3.3.7) Cordero et al. [54,55]; 3.3.8) Sharma-Arora [56]; 3.3.9) Sharma-Sharma [61]; and 3.3.10) Sharma-Guha-Gupta method [62].

Complexity of methods cannot be evaluated only based on the number of log-calls per iteration, because many of them require evaluation of derivatives in one or more points. For example, the first derivative in point $x$ is required for 2) Newton-Raphson, but not for 1) fixed-point method. On the other hand, both methods need the same number of iterations to reach the same level of accuracy. Also, in certain cases, computationally costly log-call in all subsequent iterations can be substituted with inexpensive for computation Padé polynomials [64] as described in Praks and Brkić [29], but also in that way the number of required iterations to reach the final desired accuracy remain unchanged. Therefore, for evaluation of the efficiency of each of the presented iterative methods number of iterations are counted in each of the tested 740 points within the domain of applicability and the highest value is chosen as the worst possible (five carefully selected numerical examples are chosen among 740 to illustrate the calculation).

## 3. Iterative Methods and Numerical Examples

Iterative methods with one log-call per iteration are presented in Section 3.1, with two log-calls in Section 3.2 and with three-log calls in Section 3.3. With possible exceptions, they are sorted with increased complexity. As already explained, numerical examples are given for all presented methods.

In the following formulas indexes i and i+1 refer to the two subsequent iterations.

In the following tests, the iterative procedures stop when the accuracy in respect to nine decimal places is reached or when sign #div0! appears in the meaning division with zero (desired accuracy reached). Final balanced solutions are marked with shading pattern. In framed cell are the final balanced solution that requires the highest number of iterations for the observed method (the worst case).

### *3.1. Iterative methods with one log-call per iteration*

*3.1.1*) Fixed-point method [30,31]; Equation (3)

$$x_{i+1} = x_i - F(x_i) \tag{3}$$

| | | | | |
|---|---|---|---|---|
| Example 1: | $x_1$ = 5.274011505, | $x_2$ = 5.274511624, | $x_3$ = 5.274511499, | $x_4$ = 5.274511499; |
| Example 2: | $x_1$ = 4.905054156, | $x_2$ = 4.928874894, | $x_3$ = 4.928632047, | $x_4$ = 4.928634523, |
| | $x_5$ = 4.928634490, | $x_6$ = 4.928634498, | $x_7$ = 4.928634498; | |
| Example 3: | $x_1$ = 4.128292072, | $x_2$ = 4.128359437, | $x_3$ = 4.128359435, | $x_4$ = 4.128359435; |
| Example 4: | $x_1$ = 7.331287607, | $x_2$ = 7.331277465, | $x_3$ = 7.331277467, | $x_4$ = 7.331277467; |
| Example 5: | $x_1$ = 4.124365599, | $x_2$ = 4.225356319, | $x_3$ = 4.221928724, | $x_4$ = 4.222044834, |
| | $x_5$ = 4.222040901, | $x_6$ = 4.222041034, | $x_7$ = 4.222041030, | $x_8$ = 4.222041030. |

*3.1.2*) Newton-Raphson method [32,33]; Equation (4):

$$x_{i+1} = x_i - \frac{F(x_i)}{F'(x_i)} \tag{4}$$

$F'(x_i)$ gives the fixed-point method [30,31].

| | | | | |
|---|---|---|---|---|
| Example 1: | $x_1$ = 5.274061596, | $x_2$ = 5.274511600, | $x_3$ = 5.274511499, | $x_4$ = 5.274511499; |
| Example 2: | $x_1$ = 4.907591018, | $x_2$ = 4.928826193, | $x_3$ = 4.928632752, | $x_4$ = 4.928634513, |
| | $x_5$ = 4.928634497, | $x_6$ = 4.928634498, | $x_7$ = 4.928634498; | |
| Example 3: | $x_1$ = 4.128298809, | $x_2$ = 4.128359437, | $x_3$ = 4.128359435, | $x_4$ = 4.128359435; |
| Example 4: | $x_1$ = 7.331286591, | $x_2$ = 7.331277465, | $x_3$ = 7.331277467, | $x_4$ = 7.331277467; |
| Example 5: | $x_1$ = 4.136669811, | $x_2$ = 4.224588192, | $x_3$ = 4.221965175, | $x_4$ = 4.222043289, |
| | $x_5$ = 4.222040962, | $x_6$ = 4.222041032, | $x_7$ = 4.222041030, | $x_8$ = 4.222041030. |

*3.1.3*) Hansen-Patrick method [34] is represents here through *3.1.3a*) Halley [35], *3.1.3b*) Euler-Chebyshev [36,37] and *3.1.3c*) Basto, Semiao and Calheiros methods [38]:

*3.1.3a*) Halley [35]; Equation (5):

$$x_{i+1} = x_i - \frac{\frac{F(x_i)}{F'(x_i)}}{1 - \frac{F''(x_i)}{2 \cdot F'(x_i)} \cdot \frac{F(x_i)}{F'(x_i)}} \tag{5}$$

| | | | | |
|---|---|---|---|---|
| Example 1: | $x_1$ = 5.274061858, | $x_2$ = 5.274511600, | $x_3$ = 5.274511499, | $x_4$ = 5.274511499; |
| Example 2: | $x_1$ = 4.908003000, | $x_2$ = 4.928822465, | $x_3$ = 4.928632785, | $x_4$ = 4.928634513, |
| | $x_5$ = 4.928634497, | $x_6$ = 4.928634498, | $x_7$ = 4.928634498; | |
| Example 3: | $x_1$ = 4.128298811, | $x_2$ = 4.128359437, | $x_3$ = 4.128359435, | $x_4$ = 4.128359435; |
| Example 4: | $x_1$ = 7.331286721, | $x_2$ = 7.331277465, | $x_3$ = 7.331277467, | $x_4$ = 7.331277467; |
| Example 5: | $x_1$ = 4.140502543, | $x_2$ = 4.224478344, | $x_3$ = 4.221968451, | $x_4$ = 4.222043191, |
| | $x_5$ = 4.222040965, | $x_6$ = 4.222041032, | $x_7$ = 4.222041030, | $x_8$ = 4.222041030. |

*3.1.3b*) Euler-Chebyshev method [36,37]; Equation (6):

$$x_{i+1} = x_i - \frac{F(x_i)}{F\prime(x_i)} - \frac{f^2(x_i)\cdot F''(x_i)}{2\cdot\left(F\prime(x_i)\right)^3} \tag{6}$$

Example 1: $x_1 = 5.274061740$, $x_2 = 5.274511600$, $x_3 = 5.274511499$, $x_4 = 5.274511499$;
Example 2: $x_1 = 4.907907814$, $x_2 = 4.928823333$, $x_3 = 4.928632778$, $x_4 = 4.928634513$, $x_5 = 4.928634497$, $x_6 = 4.928634498$, $x_7 = 4.928634498$;
Example 3: $x_1 = 4.128298812$, $x_2 = 4.128359437$, $x_3 = 4.128359435$, $x_4 = 4.128359435$;
Example 4: $x_1 = 7.331286591$, $x_2 = 7.331277465$, $x_3 = 7.331277467$, $x_4 = 7.331277467$;
Example 5: $x_1 = 4.141841176$, $x_2 = 4.224438148$, $x_3 = 4.221969647$, $x_4 = 4.222043156$, $x_5 = 4.222040966$, $x_6 = 4.222041032$, $x_7 = 4.222041030$, $x_8 = 4.222041030$.

*3.1.3c*) Basto-Semiao-Calheiros method [38]; Equation (7):

$$x_{i+1} = x_i - \frac{F(x_i)}{F\prime(x_i)} - \frac{f^2(x_i)\cdot F''(x_i)}{2\cdot F\prime(x_i)\cdot\left(\left(F\prime(x_i)\right)^2 - F(x_i)\cdot F''(x_i)\right)} \tag{7}$$

Example 1: $x_1 = 5.274061452$, $x_2 = 5.274511600$, $x_3 = 5.274511499$, $x_4 = 5.274511499$;
Example 2: $x_1 = 4.907358658$, $x_2 = 4.928828337$, $x_3 = 4.928632732$, $x_4 = 4.928634514$, $x_5 = 4.928634497$, $x_6 = 4.928634498$, $x_7 = 4.928634498$;
Example 3: $x_1 = 4.128298808$, $x_2 = 4.128359437$, $x_3 = 4.128359435$, $x_4 = 4.128359435$;
Example 4: $x_1 = 7.331286591$, $x_2 = 7.331277465$, $x_3 = 7.331277467$, $x_4 = 7.331277467$;
Example 5: $x_1 = 4.134234684$, $x_2 = 4.224665949$, $x_3 = 4.221962865$, $x_4 = 4.222043358$, $x_5 = 4.222040960$, $x_6 = 4.222041032$, $x_7 = 4.222041030$, $x_8 = 4.222041030$.

*3.1.4*) Super Halley method [39]; Equation (8):

$$\left.\begin{aligned} x_{i+1} &= x_i - \left(1 + \frac{1}{2}\cdot\frac{L}{1-L}\right)\cdot\frac{F(x_i)}{F\prime(x_i)} \\ L &= \frac{F(x_i)\cdot F\prime\prime(x_i)}{\left(F\prime(x_i)\right)^2} \end{aligned}\right\} \tag{8}$$

Auxiliary parameter $L$ is introduced as described.

Example 1: $x_1 = 5.274061740$, $x_2 = 5.274511600$, $x_3 = 5.274511499$, $x_4 = 5.274511499$;
Example 2: $x_1 = 4.907907729$, $x_2 = 4.928823333$, $x_3 = 4.928632778$, $x_4 = 4.928634513$, $x_5 = 4.928634497$, $x_6 = 4.928634498$, $x_7 = 4.928634498$;
Example 3: $x_1 = 4.128298812$, $x_2 = 4.128359437$, $x_3 = 4.128359435$, $x_4 = 4.128359435$, $x_5 = 4.128359435$;
Example 4: $x_1 = 7.331286591$, $x_2 = 7.331277465$, $x_3 = 7.331277467$, $x_4 = 7.331277467$;
Example 5: $x_1 = 4.141824182$, $x_2 = 4.224438659$, $x_3 = 4.221969632$, $x_4 = 4.222043156$, $x_5 = 4.222040966$, $x_6 = 4.222041032$, $x_7 = 4.222041030$, $x_8 = 4.222041030$.

*3.1.5*) Murakami [46]; Equation (9):

$$\left.\begin{aligned} x_{i+1} &= x_i - 0.3\cdot\frac{F(x_i)}{F\prime(x_i)} + \frac{1}{2}\cdot\frac{F(x_i)}{F\prime(\omega_i)} - \frac{2}{3}\cdot\frac{F(x_i)}{F\prime(\eta_i)} - \frac{32\cdot F(x_i)}{75\cdot F\prime(\omega_i) - 15\cdot F(x_i)} \\ \omega_i &= x_i - \frac{F(x_i)}{F\prime(x_i)} \\ \eta_i &= x_i - \frac{1}{2}\cdot\frac{F(x_i)}{F\prime(x_i)} \end{aligned}\right\} \tag{9}$$

In addition to the point $x$, Murakami method requires additional evaluation of the function $F$ in points $\omega$ and $\eta$ but only for the first derivative which does not contain logarithmic function.

Example 1: $x_1 = 4.918789857$, $x_2 = 5.226554297$, $x_3 = 5.269211229$, $x_4 = 5.273944807$, $x_5 = 5.274451138$, $x_6 = 5.274505072$, $x_7 = 5.274510815$, $x_8 = 5.274511426$, $x_9 = 5.274511491$, $x_{10} = 5.274511498$, $x_{11} = 5.274511499$, $x_{12} = 5.274511499$;
Example 2: $x_1 = 4.253158283$, $x_2 = 4.827210463$, $x_3 = 4.917765778$, $x_4 = 4.927553367$, $x_5 = 4.928527872$, $x_6 = 4.928623991$, $x_7 = 4.928633462$, $x_8 = 4.928634396$, $x_9 = 4.928634487$, $x_{10} = 4.928634497$, $x_{11} = 4.928634497$;

Example 3: $x_1 = 2.187900719$, $x_2 = 3.689925107$, $x_3 = 4.066519315$, $x_4 = 4.121441863$, $x_5 = 4.127617602$, $x_6 = 4.128280272$, $x_7 = 4.128350992$, $x_8 = 4.128358535$, $x_9 = 4.128359339$, $x_{10} = 4.128359425$, $x_{11} = 4.128359434$, $x_{12} = 4.128359435$, $x_{13} = 4.128359435$;

Example 4: $x_1 = 7.324855661$, $x_2 = 7.330589866$, $x_3 = 7.33120418$, $x_4 = 7.331269659$, $x_5 = 7.331276635$, $x_6 = 7.331277378$, $x_7 = 7.331277457$, $x_8 = 7.331277466$, $x_9 = 7.331277467$, $x_{10} = 7.331277467$;

Example 5: $x_1 = 2.218159942$, $x_2 = 3.806239564$, $x_3 = 4.174433434$, $x_4 = 4.21802767$, $x_5 = 4.221718266$, $x_6 = 4.222015179$, $x_7 = 4.22203896$, $x_8 = 4.222040864$, $x_9 = 4.222041017$, $x_{10} = 4.222041029$, $x_{11} = 4.222041030$, $x_{12} = 4.222041030$.

*3.2. Iterative methods with two log-calls per iteration*

*3.2.1)* Ostrowski method of King family [36,37,40]; Equation (10):

$$\left.\begin{aligned} x_{i+1} &= y_i - \frac{F(y_i)}{F'(x_i)} \cdot \frac{F(x_i)}{F(x_i) - 2 \cdot F(y_i)} \\ y_i &= x_i - \frac{F(x_i)}{F'(x_i)} \end{aligned}\right\} \quad (10)$$

Example 1: $x_1 = 5.274511398$, $x_2 = 5.274511499$, $x_3 =$ #div0;
Example 2: $x_1 = 4.928451807$, $x_2 = 4.928634512$, $x_3 = 4.928634498$, $x_4 = 4.928634498$;
Example 3: $x_1 = 4.128359434$, $x_2 = 4.128359435$, $x_3 =$ #div0;
Example 4: $x_1 = 7.331277468$, $x_2 = 7.331277467$, $x_3 =$ #div0;
Example 5: $x_1 = 4.219926077$, $x_2 = 4.222042800$, $x_3 = 4.222041028$, $x_4 = 4.222041030$, $x_5 = 4.222041030$.

*3.2.2)* Kung-Traub method [41]; Equation (11):

$$\left.\begin{aligned} x_{i+1} &= y_i - \frac{F(y_i)}{F'(x_i)} \cdot \frac{1}{\left(1 - \frac{F(y_i)}{F(x_i)}\right)^2} \\ y_i &= x_i - \frac{F(x_i)}{F'(x_i)} \end{aligned}\right\} \quad (11)$$

Example 1: $x_1 = 5.274511398$, $x_2 = 5.274511499$, $x_3 =$ #div0!;
Example 2: $x_1 = 4.928450156$, $x_2 = 4.928634513$, $x_3 = 4.928634498$, $x_4 = 4.928634498$;
Example 3: $x_1 = 4.128359434$, $x_2 = 4.128359435$, $x_3 =$ #div0!;
Example 4: $x_1 = 7.331277468$, $x_2 = 7.331277467$, $x_3 =$ #div0!;
Example 5: $x_1 = 4.219864191$, $x_2 = 4.222042905$, $x_3 = 4.222041028$, $x_4 = 4.222041030$, $x_5 = 4.222041030$.

*3.2.3)* Maheshwari method [42]; Equation (12):

$$\left.\begin{aligned} x_{i+1} &= x_i - \left(\left(\frac{F(y_i)}{F(x_i)}\right)^2 - \frac{F(x_i)}{F(y_i) - F(x_i)}\right) \cdot \frac{F(x_i)}{F'(x_i)} \\ y_i &= x_i - \frac{F(x_i)}{F'(x_i)} \end{aligned}\right\} \quad (12)$$

Example 1: $x_1 = 5.274511398$, $x_2 = 5.274511499$, $x_3 =$ #div0!;
Example 2: $x_1 = 4.928446781$, $x_2 = 4.928634513$, $x_3 = 4.928634498$, $x_4 = 4.928634498$;
Example 3: $x_1 = 4.128359434$, $x_2 = 4.128359435$, $x_3 =$ #div0!;
Example 4: $x_1 = 7.331277468$, $x_2 = 7.331277467$, $x_3 =$ #div0!;
Example 5: $x_1 = 4.219731647$, $x_2 = 4.222043139$, $x_3 = 4.222041028$, $x_4 = 4.222041030$, $x_5 = 4.222041030$.

*3.2.4)* Khattri and Babajee [45]; Equation (13):

$$\left.\begin{aligned} x_{i+1} &= y_i - \frac{F(x_i) \cdot F(y_i)}{F(x_i) - 2 \cdot F(y_i)} \cdot \left(\frac{3}{F'(x_i) + 0.001 \cdot F(y_i)} - \frac{2}{F'(x_i)}\right) \\ y_i &= x_i - \frac{F(x_i)}{F'(x_i)} \end{aligned}\right\} \quad (13)$$

Example 1: $x_1 = 5.274511397$, $x_2 = 5.274511499$, $x_3 = \#div0!$;
Example 2: $x_1 = 4.928450478$, $x_2 = 4.928634513$, $x_3 = 4.928634498$, $x_4 = 4.928634498$;
Example 3: $x_1 = 4.128359434$, $x_2 = 4.128359435$, $x_3 = \#div0!$;
Example 4: $x_1 = 7.331277468$, $x_2 = 7.331277467$, $x_3 = \#div0!$;
Example 5: $x_1 = 4.219904119$, $x_2 = 4.222042819$, $x_3 = 4.222041028$, $x_4 = 4.222041030$, $x_5 = 4.222041030$.

*3.2.5*) Hermite interpolation [43] based on Jarratt method [44]; Equation (14):

In Equation (14), $H$ is Hermite interpolation polynomial. This method requires evaluation of the first derivatives in points $x$, $y$ and $z$, but the function $F$ need to be evaluated only in points $x$ and $y$.

$$\left.\begin{aligned} x_{i+1} &= z_i - \frac{H_i}{F'(z_i)} \\ y_i &= x_i - \frac{2}{3}\cdot\frac{F(x_i)}{F'(x_i)} \\ z_i &= x_i - \frac{1}{2}\cdot\frac{F(x_i)}{F'(x_i)}\cdot\left(1+\frac{1}{1+\frac{3}{2}\cdot\left(\frac{F'(y_i)}{F'(x_i)}-1\right)}\right) \\ H_i &= F(x_i) + F'(x_i)\cdot\frac{(z_i-x_i)\cdot(z_i-y_i)^2}{(y_i-x_i)\cdot(x_i+2\cdot y_i-3\cdot z_i)} + F'(z_i)\cdot\frac{(z_i-y_i)\cdot(x_i-z_i)}{x_i+2\cdot y_i-3\cdot z_i} - \frac{F(x_i)-F(y_i)}{x_i-y_i}\cdot\frac{(z_i-x_i)^3}{(y_i-x_i)\cdot(x_i+2\cdot y_i-3\cdot z_i)} \end{aligned}\right\} \quad (14)$$

Example 1: $x_1 = 5.274466557$, $x_2 = 5.2745115$, $x_3 = 5.274511499$, $x_4 = 5.274511499$;
Example 2: $x_1 = 4.926606155$, $x_2 = 4.928636193$, $x_3 = 4.928634496$, $x_4 = 4.928634498$, $x_5 = 4.928634498$;
Example 3: $x_1 = 4.128353373$, $x_2 = 4.128359436$, $x_3 = 4.128359435$, $x_4 = \#div0!$;
Example 4: $x_1 = 7.331278378$, $x_2 = 7.331277467$, $x_3 = 7.331277467$;
Example 5: $x_1 = 4.214067429$, $x_2 = 4.222058401$, $x_3 = 4.222040992$, $x_4 = 4.222041030$, $x_5 = 4.222041030$.

*3.2.6*) Wang-Liu method [49] based on Hermit interpolation [43]; Equation (15):

In Equation (15), $H$ is Hermite interpolation polynomial used in the same form as in Equation (14). This method also requires evaluation of the first derivatives in point $x$, $y$ and $z$, but the function $F$ need to be evaluated only in points $x$ and $y$.

$$\left.\begin{aligned} x_{i+1} &= z_i - \frac{H_i}{F'(z_i)} \\ y_i &= x_i - \frac{F(x_i)}{F'(x_i)} \\ z_i &= y_i - \frac{F(y_i)}{F'(x_i)}\cdot\frac{F(x_i)}{F(x_i)-2\cdot F(y_i)} \\ H_i &= F(x_i) + F'(x_i)\cdot\frac{(z_i-x_i)\cdot(z_i-y_i)^2}{(y_i-x_i)\cdot(x_i+2\cdot y_i-3\cdot z_i)} + F'(z_i)\cdot\frac{(z_i-y_i)\cdot(x_i-z_i)}{x_i+2\cdot y_i-3\cdot z_i} - \frac{F(x_i)-F(y_i)}{x_i-y_i}\cdot\frac{(z_i-x_i)^3}{(y_i-x_i)\cdot(x_i+2\cdot y_i-3\cdot z_i)} \end{aligned}\right\} \quad (15)$$

Example 1: $x_1 = 5.274061596$, $x_2 = 5.2745116$, $x_3 = 5.274511499$, $x_4 = 5.274511499$;
Example 2: $x_1 = 4.907590959$, $x_2 = 4.928826193$, $x_3 = 4.928632752$, $x_4 = 4.928634513$, $x_5 = 4.928634497$, $x_6 = 4.928634498$, $x_7 = 4.928634498$;
Example 3: $x_1 = 4.128298809$, $x_2 = 4.128359437$, $x_3 = 4.128359435$, $x_4 = 4.128359435$;
Example 4: $x_1 = 7.331286591$, $x_2 = 7.331277465$, $x_3 = 7.331277467$, $x_4 = 7.331277467$;
Example 5: $x_1 = 4.136665933$, $x_2 = 4.224588276$, $x_3 = 4.221965174$, $x_4 = 4.222043289$, $x_5 = 4.222040962$, $x_6 = 4.222041032$, $x_7 = 4.222041030$, $x_8 = 4.222041030$.

### *3.3. Iterative methods with three log-calls per iteration*

*3.3.1*) Neta [47]; Equation (16):

$$\left.\begin{aligned} x_{i+1} &= z_i - \frac{F(z_i)}{F'(x_i)}\cdot\frac{F(x_i)-F(y_i)}{F(x_i)-3\cdot F(y_i)} \\ y_i &= x_i - \frac{F(x_i)}{F'(x_i)} \\ z_i &= y_i - \frac{F(y_i)}{F'(x_i)}\cdot\frac{F(x_i)-\frac{1}{2}\cdot F(y_i)}{F(x_i)-\frac{5}{2}\cdot F(y_i)} \end{aligned}\right\} \quad (16)$$

Example 1: $x_1 = 5.274511499$, $x_2 = 5.274511499$;
Example 2: $x_1 = 4.928632954$, $x_2 = 4.928634498$, $x_3 = 4.928634498$;
Example 3: $x_1 = 4.128359435$, $x_2 = 4.128359435$;
Example 4: $x_1 = 7.331277467$, $x_2 = 7.331277467$;
Example 5: $x_1 = 4.221992945$, $x_2 = 4.222041031$, $x_3 = 4.222041031$.

*3.3.2*) Chun-Neta [48] based on Kung-Traub [41]; Equation (17):

$$\left.\begin{aligned} x_{i+1} &= z_i - \frac{F(z_i)}{F\prime(x_i)} \cdot \frac{1}{\left(1 - \frac{F(y_i)}{F(x_i)} - \frac{F(z_i)}{F(x_i)}\right)^2} \\ y_i &= x_i - \frac{F(x_i)}{F\prime(x_i)} \\ z_i &= y_i - \frac{F(y_i)}{F\prime(x_i)} \cdot \frac{1}{\left(1 - \frac{F(y_i)}{F(x_i)}\right)^2} \end{aligned}\right\} \tag{17}$$

Example 1: $x_1 = 5.274511499$, $x_2 = 5.274511499$;
Example 2: $x_1 = 4.928632854$, $x_2 = 4.928634498$, $x_3 = 4.928634498$;
Example 3: $x_1 = 4.128359435$, $x_2 = 4.128359435$;
Example 4: $x_1 = 7.331277467$, $x_2 = 7.331277467$;
Example 5: $x_1 = 4.221982464$, $x_2 = 4.222041031$, $x_3 = 4.222041030$, $x_4 = 4.222041030$.

*3.3.3*) Džunić-Petković-Petković [57-59]; Equation (18):

$$\left.\begin{aligned} x_{i+1} &= z_i - \frac{F(z_i)}{F'(x_i) \cdot \left[1 - 2 \cdot \frac{F(y_i)}{F(x_i)} - \left(\frac{F(y_i)}{F(x_i)}\right)^2\right] \cdot \left[1 - \frac{F(z_i)}{F(y_i)}\right] \cdot \left[1 - 2 \cdot \frac{F(z_i)}{F(x_i)}\right]} \\ y_i &= x_i - \frac{F(x_i)}{F\prime(x_i)} \\ z_i &= y_i - \frac{F(x_i)}{F(x_i) - 2 \cdot F(y_i)} \cdot \frac{F(y_i)}{F'(x_i)} \end{aligned}\right\} \tag{18}$$

Example 1: $x_1 = 5.274511499$, $x_2 =$ #div0;
Example 2: $x_1 = 4.928634483$, $x_2 = 4.928634498$, $x_3 =$ #div0;
Example 3: $x_1 = 4.128359435$, $x_2 =$ #div0;
Example 4: $x_1 = 7.331277467$, $x_2 =$ #div0;
Example 5: $x_1 = 4.222039554$, $x_2 = 4.222041030$, $x_3 = 4.222041030$.

*3.3.4*) Neta-Johnson [60] based on Jarrat method [44]; Equation (19):

$$\left.\begin{aligned} x_{i+1} &= z_n - \frac{F(z_i)}{F'(x_i)} \cdot \frac{F'(x_i) + F'(y_i) - F'(\delta_i)}{-2 \cdot F'(x_i) + 2 \cdot F'(y_i) - F'(\delta_i)} \\ y_i &= x_i - \frac{F(x_i)}{F\prime(x_i)} \\ \delta_i &= x_i - \frac{1}{8} \cdot \frac{F(x_i)}{F'(x_i)} - \frac{3}{8} \cdot \frac{F(x_i)}{F(y_i)} \\ z_i &= x_i - \frac{F(x_i)}{\frac{1}{6} \cdot F'(x_i) + \frac{1}{6} \cdot F'(y_i) + \frac{2}{3} \cdot F'(\delta_i)} \end{aligned}\right\} \tag{19}$$

Auxiliary parameter $\delta$ is introduced as described.

Example 1: $x_1 = 5.272711608$, $x_2 = 5.27451312$, $x_3 = 5.274511498$, $x_4 = 5.274511499$, $x_5 = 5.274511499$;
Example 2: $x_1 = 4.843943256$, $x_2 = 4.931741639$, $x_3 = 4.928520569$, $x_4 = 4.928638675$, $x_5 =4.928634344$, $x_6 = 4.928634503$, $x_7 = 4.928634497$, $x_8 = 4.928634498$, $x_9 = 4.928634498$;
Example 3: $x_1 = 4.128116927$, $x_2 = 4.128359454$, $x_3 = 4.128359435$, $x_4 = 4.128359435$;
Example 4: $x_1 = 7.331313968$, $x_2 = 7.331277444$, $x_3 = 7.331277467$, $x_4 = 7.331277467$;
Example 5: $x_1 = 3.873979004$, $x_2 = 4.264698085$, $x_3 = 4.21685805$, $x_4 = 4.222671442$, $x_5 = 4.221964362$, $x_6 =4.222050354$, $x_7 = 4.222039896$, $x_8 = 4.222041168$, $x_9 =4.222041013$, $x_{10} = 4.222041032$, $x_{11} =4.222041030$, $x_{12} = 4.222041030$.

*3.3.5*) Jain [51] based on Steffensen method [52-55]; Equation (20):

$$\left.\begin{aligned} x_{i+1} &= x_i - \frac{F^3(x_i)}{\big(F(x_i+F(x_i))-F(x_i)\big)\cdot\big(F(x_i)-F(y_i)\big)} \\ y_i &= x_i - \frac{F^2(x_i)}{F(x_i+F(x_i))-F(x_i)} \end{aligned}\right\} \tag{20}$$

Jain method is three-point method with evaluation of the function $F$ in points $x$, $y$, but also in point $x_i + F(x_i)$. Serghides explicit approximation of the Colebrook equation is based on Steffensen method [16, 65].

Example 1: $x_1 = 5.274511499$, $x_2 = 5.274511499$;
Example 2: $x_1 = 4.928634582$, $x_2 = 4.928634498$, $x_3 =$ #div0;
Example 3: $x_1 = 4.128359435$, $x_2 =$ #div0;
Example 4: $x_1 = 7.331277467$, $x_2 =$ #div0;
Example 5: $x_1 = 4.222058673$, $x_2 = 4.222041030$, $x_3 =$ #div0.

*3.3.6*) Bi-Ren-Wu [50]; Equation (21):

$$\left.\begin{aligned} x_{i+1} &= z_i - \frac{F(z_i)}{(\lozenge_{zy})_i+(\lozenge_{yx})_i-F'(x_i)} \\ y_i &= x_i - \frac{F(x_i)}{F'(x_i)} \\ z_i &= y_i - \frac{F(y_i)}{F'(x_i)}\cdot\frac{F(x_i)}{F(x_i)-2\cdot F(y_i)} \\ (\lozenge_{yx})_i &= \frac{F(y_i)-F(x_i)}{y_i-x_i} \\ (\lozenge_{zy})_i &= \frac{F(z_i)-F(y_i)}{z_i-y_i} \end{aligned}\right\} \tag{21}$$

Two-parameter function $\lozenge$ is introduced as defined.

Example 1: $x_1 = 5.274511432$, $x_2 = 5.274511499$, $x_3 =$ #div0!;
Example 2: $x_1 = 4.928634504$, $x_2 = 4.928634498$, $x_3 =$ #div0!;
Example 3: $x_1 = 4.128359435$, $x_2 = 4.128359435$;
Example 4: $x_1 = 7.331277468$, $x_2 = 7.331277467$, $x_3 =$ #div0!;
Example 5: $x_1 = 4.222041809$, $x_2 = 4.222041029$, $x_3 = 4.222041030$, $x_4 = 4.222041030$.

*3.3.7*) Cordero et al. [54,55]; Equation (22):

$$\left.\begin{aligned} x_{i+1} &= z_i - \frac{1+3\cdot\frac{F(z_i)}{F(x_i)}}{1+\frac{F(z_i)}{F(x_i)}}\cdot\frac{F(z_i)}{(\lozenge_{zy})_i+(\lozenge_{zxx})_i\cdot(z_i-y_i)} \\ y_i &= x_i - \frac{F(x_i)}{F'(x_i)} \\ z_i &= y_i - \frac{F(y_i)}{F'(x_i)}\cdot\frac{1}{1-2\cdot\frac{F(y_i)}{F(x_i)}-\left(\frac{F(y_i)}{F(x_i)}\right)^2-\frac{\left(\frac{F(y_i)}{F(x_i)}\right)^3}{2}} \\ (\lozenge_{zy})_i &= \frac{F(z_i)-F(y_i)}{z_i-y_i} \\ (\lozenge_{zx})_i &= \frac{F(z_i)-F(x_i)}{z_i-x_i} \\ (\lozenge_{zxx})_i &= \frac{(\lozenge_{zx})_i-F'(x_i)}{z_i-x_i} \end{aligned}\right\} \tag{22}$$

Two-parameter and three-parameter functions $\lozenge$ are introduced as defined.

Example 1: $x_1 = 5.274511398$, $x_2 = 5.274511499$, $x_3 =$ #div0!;
Example 2: $x_1 = 4.928453453$, $x_2 = 4.928634512$, $x_3 = 4.928634498$, $x_4 = 4.928634498$;
Example 3: $x_1 = 4.128359434$, $x_2 = 4.128359435$, $x_3 =$ #div0!;
Example 4: $x_1 = 7.331277468$, $x_2 = 7.331277467$, $x_3 =$ #div0!;
Example 5: $x_1 = 4.219987477$, $x_2 = 4.2220427$, $x_3 = 4.222041028$, $x_4 = 4.222041030$, $x_5 = 4.222041030$.

*3.3.8*) Sharma-Arora [56]; Equation (23):

$$\left.\begin{aligned} x_{i+1} &= z_i - \frac{(\lozenge_{zy})_i}{(\lozenge_{zx})_i} \cdot \frac{F(z_i)}{2\cdot(\lozenge_{zy})_i - (\lozenge_{zx})_i} \\ y_i &= x_i - \frac{F(x_i)}{F\prime(x_i)} \\ (\lozenge_{yx})_i &= \frac{F(y_i)-F(x_i)}{y_i - x_i} \\ z_i &= y_i - \frac{F(y_i)}{2\cdot(\lozenge_{yx})_i - F\prime(x_i)} \\ (\lozenge_{zx})_i &= \frac{F(z_i)-F(x_i)}{z_i - x_i} \\ (\lozenge_{zy})_i &= \frac{F(z_i)-F(y_i)}{z_i - y_i} \end{aligned}\right\} \quad (23)$$

Two-parameter function $\lozenge$ is introduced as defined.

| | | | |
|---|---|---|---|
| Example 1: | $x_1 = 5.274511499,$ | $x_2 = \#div0!;$ | |
| Example 2: | $x_1 = 4.928634497,$ | $x_2 = 4.928634498,$ | $x_3 = \#div0!;$ |
| Example 3: | $x_1 = 4.128359435,$ | $x_2 = \#div0!;$ | |
| Example 4: | $x_1 = 7.331277467,$ | $x_2 = \#div0!;$ | |
| Example 5: | $x_1 = 4.222040921,$ | $x_2 = 4.222041030,$ | $x_3 = \#div0!.$ |

*3.3.9*) Sharma-Sharma [61]; Equation (24):

$$\left.\begin{aligned} x_{i+1} &= z_i - w_i \cdot \frac{F(z_i)\cdot(\lozenge_{xy})_i}{(\lozenge_{xz})_i\cdot(\lozenge_{yz})_i} \\ y_i &= x_i - \frac{F(x_i)}{F\prime(x_i)} \\ z_i &= y_i - \frac{F(y_i)}{F'(x_i)} \cdot \frac{1}{1-2\cdot\frac{F(y_i)}{F(x_i)}} \\ w_i &= 1 + \frac{\frac{F(z_i)}{F(x_i)}}{1+\frac{F(z_i)}{F(x_i)}} \\ (\lozenge_{xy})_i &= \frac{F(x_i)-F(y_i)}{x_i - y_i} \\ (\lozenge_{xz})_i &= \frac{F(x_i)-F(z_i)}{x_i - z_i} \\ (\lozenge_{yz})_i &= \frac{F(y_i)-F(z_i)}{y_i - z_i} \end{aligned}\right\} \quad (24)$$

Two-parameter function $\lozenge$ is introduced as defined.

| | | | |
|---|---|---|---|
| Example 1: | $x_1 = 5.274511499,$ | $x_2 = \#div0!;$ | |
| Example 2: | $x_1 = 4.928634483,$ | $x_2 = 4.928634498,$ | $x_3 = \#div0!;$ |
| Example 3: | $x_1 = 4.128359435,$ | $x_2 = \#div0!;$ | |
| Example 4: | $x_1 = 7.331277467,$ | $x_2 = \#div0!;$ | |
| Example 5: | $x_1 = 4.222039549,$ | $x_2 = 4.222041030,$ | $x_3 = 4.22204103.$ |

*3.3.10*) Sharma-Guha-Gupta method [62]; Equation (25):

$$\left.\begin{aligned} x_{i+1} &= x_i - \frac{P+Q+R}{P\cdot(\lozenge_{zx})_i + Q\cdot F'(x_i) + R\cdot(\lozenge_{yx})_i} \cdot F(x_i) \\ y_i &= x_i - \frac{F(x_i)}{F\prime(x_i)} \\ z_i &= y_i - \frac{1}{1-2\cdot\frac{F(y_i)}{F(x_i)}} \cdot \frac{F(y_i)}{F'(x_i)} \\ P &= (x_i - y_i)\cdot F(x_i)\cdot F(y_i) \\ Q &= (y_i - z_i)\cdot F(y_i)\cdot F(z_i) \\ R &= (z_i - x_i)\cdot F(z_i)\cdot F(x_i) \\ (\lozenge_{zx})_i &= \frac{F(z_i)-F(x_i)}{z_i - x_i} \\ (\lozenge_{yx})_i &= \frac{F(y_i)-F(x_i)}{y_i - x_i} \end{aligned}\right\} \quad (25)$$

Two-parameter function ◊ is introduced as defined, as well as auxiliary parameters $P$, $Q$ and $R$.

Example 1: $x_1 = 5.274511499,$ $x_2 = \#div0!;$
Example 2: $x_1 = 4.928634483,$ $x_2 = 4.928634498,$ $x_3 = \#div0!;$
Example 3: $x_1 = 4.128359435,$ $x_2 = \#div0!;$
Example 4: $x_1 = 7.331277467,$ $x_2 = \#div0!;$
Example 5: $x_1 = 4.222039558,$ $x_2 = 4.222041030,$ $x_3 = 4.222041030.$

## 4. Summary - Discussion and Analysis

To summarize, the highest required number of iterations to reach the final balanced solution of the Colebrook equation in respect to nine decimal places, for the examined methods are:

**3.1 One log-call per iteration (one-point methods)**

3.1.1) Fixed-point [30,31]; Equation (3): **7 iterations**
3.1.2) Newton-Raphson [32,33]; Equation (4): **7 iterations**
3.1.3) Hansen-Patrick [34]:
- 3.1.3a) Halley [35], Equation (5): **7 iterations**
- 3.1.3b) Euler-Chebyshev [36,37], Equation (6): **7 iterations**
- 3.1.3c) Basto-Semiao-Calheiros method [38]; Equation (7): **7 iterations**

3.1.4) Super Halley [39]; Equation (8): **7 iterations**
3.1.5) Murakami [46]; Equation (9): **12 iterations**

**3.2 Two log-calls per iteration (two-point methods):**

3.2.1) Ostrowski method of King family [36,37,40]; Equation (10): **4 iterations**
3.2.2) Kung-Traub [41]; Equation (11): **4 iterations**
3.2.3) Maheshwari [42]; Equation (12): **4 iterations**
3.2.4) Khattri and Babajee [45]; Equation (13): **4 iterations**
3.2.5) Hermite interpolation [43] based on Jarratt method [44]; Equation (14): **4 iterations**
3.2.6) Wang-Liu method [49] based on Hermit interpolation [43]; Equation (15): **7 iterations**

**3.3 Three log-calls per iteration (three-point methods):**

3.3.1) Neta [47]; Equation (16): **2 iterations**
3.3.2) Chun-Neta [48] based on Kung-Traub [41]; Equation (17): **2 iterations** (3 in rare cases)
3.3.3) Džunić-Petković-Petković method [57-59]; Equation (18): **2 iterations**
3.3.4) Neta-Johnson [60] based on Jarrat method [44]; Equation (19): **11 iterations**
3.3.5) Jain [51] based on Steffensen method [52-55]; Equation (20): **2 iterations**
3.3.6) Bi-Ren-Wu [50]; Equation (21): **3 iterations**
3.3.7) Cordero et al. [54,55]; Equation (22): **4 iterations**
3.3.8) Sharma-Arora [56]; Equation (23): **2 iterations**
3.3.9) Sharma-Sharma [61]; Equation (24): **2 iterations**
3.3.10) Sharma-Guha-Gupta method [62]; Equation (25): **2 iterations**

After conducted analysis the following methods should not be used for the Colebrook equation: 3.1.5) Murakami [46]; Equation (9), 3.2.6) Wang-Liu method [49] based on Hermit interpolation [43]; Equation (15), Neta-Johnson [60] based on Jarrat method [44]; Equation (19), and Cordero et al. [54,55]; Equation (22).

Most one-point methods that require one log-call per iteration need seven iterations to reach final solution. Among them as the most simplest, 3.1.1) Fixed-point [30,31]; Equation (3) can be recommended. Following procedure from Praks and Brkić [29], only one log-call is required in respect to the whole procedure and in particular only in the first iteration, where in all subsequent iterations instead of log-call as substitution can be used computationally inexpensive Padé polynomials.

Two log-calls iterative methods require up to four iterations to reach final solution. Those simpler can be used: 3.2.1) Ostrowski method of King family [36,37,40]; Equation (10), 3.2.2) Kung-Traub [41]; Equation (11), and 3.2.3) Maheshwari [42]; Equation (12).

Three log-calls methods are the most powerful among the presented procedures. The simplest, but accurate, can be recommended for use: 3.3.1) Neta [47]; Equation (16), 3.3.2) Chun-Neta [48] based on Kung-Traub [41]; Equation (17), 3.3.3) Džunić-Petković-Petković method [57-59]; Equation (18), and 3.3.5) Jain [51] based on Steffensen method [52-55]; Equation (20). Those accurate, but with two-parameter function ◊ can be used, but for the Colebrook equation they seem to be too complex and their use should be limited to some rare very well elaborated cases and calculations: Sharma-Arora [56]; Equation (23), 3.3.9) Sharma-Sharma [61]; Equation (24), and 3.3.10) Sharma-Guha-Gupta method [62]; Equation (25).

Choosing the different starting point $x_0$ compared with the here selected for the shown numerical validation does not alter the number of required iterations significantly [29,66].

## 5. Conclusions

The Colebrook equation for flow friction in its domain of applicability is fast converging. The fixed-point iterative methods are in common use but they demand up to seven iterations to reach the final satisfied and balanced accuracy [30,31]. On the other hand, numerous explicit approximations with different degree of accuracy are available [12-14]. Similar as here presented two-point and three-point iterative procedures, the very accurate approximations of the Colebrook equation also contain few internal iterative steps [15-17]. Here tested and proposed three-point iterative methods are very accurate, and therefore they can be used in a hybrid way: (1) as very accurate explicit approximations if the calculation is terminated after the first iteration, or (2) for very fine and punctual computing where the very high accuracy can be easily reached after the one or two more additional iterations. This is of great significance in computing of flow through miscellaneous pipe-networks of various applications (water, natural gas, air, etc.) where multiple evaluation of friction factor takes place and when certain cases requires only rough estimation, while another very accurate scientific computing that need to be compared, repeated or analyzed in details [67-69].

**Notations**

The following symbols are used in this paper:

$\lambda$-Darcy friction factor (Moody, Darcy-Weisbach or Colebrook); dimensionless

$Re$-Reynolds number; dimensionless

$\varepsilon^*$-relative roughness of inner pipe surface; dimensionless

$x=\frac{1}{\sqrt{\lambda}}$-linearization of friction factor – first point; dimensionless

$y=\frac{1}{\sqrt{\lambda}}$-linearization of friction factor – second point; dimensionless

$z=\frac{1}{\sqrt{\lambda}}$-linearization of friction factor – third point; dimensionless

$i$-counter

Auxiliary parameters

$L$- used in Super Halley method [39]; Equation (8)

$\omega$-used in Murakami [46]; Equation (9)

$\eta$- used in Murakami [46]; Equation (9)

$H$-Hermite interpolation polynomial used in Jarratt method [44]; Equation (14) and in Wang-Liu method [49]; Equation (15)

$\delta$-used in Neta-Johnson [60] based on Jarrat method [44]; Equation (19)

$P$, $Q$ and $R$-Sharma-Guha-Gupta method [62]; Equation (25)

◊-Two-parameter or three-parameter function used in Bi-Ren-Wu [50]; Equation (21), Cordero et al. [54,55]; Equation (22), Sharma-Arora [56]; Equation (23), Sharma-Sharma [61]; Equation (24), Sharma-Guha-Gupta method [62]; Equation (25)

Functional symbols:

$\zeta$-Colebrook equation

F-Colebrook equation in form; $F(x) = x - \zeta(x)$

$F^{'}$-first derivative

$F^{''}$-second derivative

$log_{10}$-Briggs logarithm

$ln$-Napier natural logarithm

**Author Contributions:** Both authors contributed equally to this study. The paper is product of joint efforts of the authors who worked together on models of natural gas distribution networks. Pavel Praks has scientific background in applied mathematics and programming while Dejan Brkić in control and applied computing in mechanical and petroleum engineering. Following the idea by Pavel Praks, Dejan Brkić tested the methods and controlled the results. Dejan Brkić wrote the draft of this paper.

**Conflicts of Interest:** The authors declare no conflict of interest. The views expressed are those of the authors and may not in any circumstances be regarded as stating an official position of the European Commission or of the Technical University Ostrava.